# Optimum dimensional synthesis of planar mechanisms with geometric constraints

**V. García-Marina** · **I. Fernández de Bustos** · **G. Urkullu** · **R. Ansola**

**Abstract** The deformed energy method has shown to be a good option for dimensional synthesis of mechanisms. In this paper the introduction of some new features to such approach is proposed. First, constraints fixing dimensions of certain links are introduced in the error function of the synthesis problem. Second, requirements on distances between determinate nodes are included in the error function for the analysis of the deformed position problem. Both the overall synthesis error function and the inner analysis error function are optimized using a Sequential Quadratic Problem (SQP) approach. This also reduces the probability of branch or circuit defects. In the case of the inner function analytical derivatives are used, while in the synthesis optimization approximate derivatives have been introduced. Furthermore, constraints are analyzed under two formulations, the Euclidean distance and an alternative approach that uses the previous raised to the power of two. The latter approach is often used in kinematics, and simplifies the computation of derivatives. Some examples are provided to show the convergence order of the error function and the fulfilment of the constraints in both formulations studied under different topological situations or achieved energy levels.

V. García-Marina (✉)
Department of Mechanical Engineering, Faculty of Engineering of Vitoria-Gasteiz, University of the Basque Country (UPV/EHU), Nieves Cano 12, 01006 Vitoria-Gasteiz, Spain
e-mail: vanessa.garcia@ehu.eus

I. Fernández de Bustos · G. Urkullu · R. Ansola
Department of Mechanical Engineering, Faculty of Engineering of Bilbao, University of the Basque Country (UPV/EHU), Alameda de Urquijo s/n, 48013 Bilbao, Spain
e-mail: igor.fernandezdebustos@ehu.es

G. Urkullu
e-mail: gorka.urkullu@ehu.eus

R. Ansola
e-mail: ruben.ansola@ehu.eus

## 1 Introduction

### 1.1 Background

The design of mechanisms is a complex field in mechanical engineering, and many methods for general mechanism design have been proposed in the literature. Some among them are found in [1] for kinematic synthesis, an Assur Group based method for type kinematic synthesis and Genetic Algorithms for dimensional kinematic synthesis as in [2], an Exact

Gradient method in [3] for dimensional kinematic synthesis, Neural Networks for dimensional kinematic synthesis [4], Geometric Constraint Programming for kinematic synthesis as in [5] and [6] among others, or Kinematic Mapping for kinematic synthesis in [7], Circular Proximity Function for optimal synthesis in [8], or synthesis of constrained chains in [9]. Given that it is undesirable to work with a different program for each type of mechanism or problem, attempts have been made to develop a versatile method for all possible topologies of mechanisms, problems of kinematic analysis of position (initial position, successive positions, position of static equilibrium, or deformed position) such as [10, 11], or problems of synthesis (path generation, generation of functions, and rigid solid guidance) to give an example [8, 12–17].

The synthesis of mechanisms is a heavily nonlinear problem, and this usually leads to resort to iterative procedures whose convergence may not always be guaranteed, but only for certain starting conditions.

Methods to optimize mechanisms can be divided into: *Heuristics*, algorithms based on analogies to natural or artificial phenomena that have been widely used in mechanism synthesis. They include genetic algorithms [2, 18, 19], simulated annealing [20, 21], ant colony optimization [22], particle swarm optimization [23, 24], or tabu search [25] among others.

*Analytical*, not easily changeable methods, i.e., intended for the analysis of concrete and not very complex mechanisms like four or six links. Usually defined for a certain problem, and not applicable to other problems [13, 14, 26, 27].

*Numerical*, defined to solve general and complex problems with many links. Several methods are available to tackle the synthesis of mechanisms, such as the Dogleg method [28], the simplex method [29], the interior point method [30], penalty functions [16, 31], and gradient-based methods [15] [3] [32, 33]. However, variants of the SQP are used more often for mechanism synthesis [34–36].

As this study centers on SQP methods, the authors find it useful to describe the more usual variations used to implement constraints to the problem.

1. For linear equality constraints: ad-hoc methods, which systematically eliminate rows and columns, and the Null Space method [37–39].

2. For nonlinear equality constraints: Penalty Function methods [40, 41], Lagrange Multiplier methods [42, 43], and Augmented Lagrangian Function methods [44, 45].
3. For linear inequality constraints: the Karmarkar methods [46, 47] and Primal-Dual methods [48, 49].
4. For nonlinear inequality constraints: the Logarithmic Barrier Function methods [50, 51] and the Slack Variables methods [52–54].

### 1.2 Formulation of the problem of interest for this investigation

In this paper, some improvements are applied to an existing method to solve the deformed position problem, that enable it to yield an optimal dimensional synthesis and offer a kinematic analysis of a given mechanism. The implementation of the aforementioned geometric constraints in algorithms has several applications, e.g., the need for fixing the length of a certain element of a mechanism when designing machines (dimensional constraints), or when modeling such diverse elements as grippers, springs, or dampers, where the necessity to introduce requirements can be formulated in the form of a distance between two points (distance constraints). These applications show the need of such restrictions which cannot be implemented to the optimization problem with the formulations usually available for path generation, function generation or solid rigid guidance. In consequence, the proposal of these geometric constraints opens a wider field of employment for the dimensional synthesis.

Although the authors find type synthesis an interesting field to work on, they have centered this work on dimensional synthesis. Therefore, in this paper topology of the mechanism is preset.

### 1.3 Scope and contribution of this study

The method implemented in this paper, as well as the original method, is of optimization character, therefore the solution will not always fulfil exactly all the requirements, and an agreement solution will be searched for. The requirements are not always achieved with accuracy, but rather a solution is sought that is as close as possible to all points. Therefore, the solution chosen verifies the constraints as best as possible so the deformation energy stored in the mechanism will not be zero, but approximates to it.

The evaluation function used is the deformation energy stored in the mechanism when it is forced to occupy positions that do not fall within the natural range of its path without deforming its elements [55–58], assuming that the mechanism is composed of elastic and deformable elements. The method was first proposed by R. Avilés in 1982, [55], and has undergone continual advances [56, 57, 59, 60], and [61]. This function is a clear example of nonlinearity. Besides, to constitute the augmented error function, a series of geometric constraints are incorporated using the Lagrange multipliers, which is also of nonlinear character.

A significant improvement recently introduced in [62] is also used here, and implies a reconfiguration of the parameters of the synthesis. In [62] the authors presented the method based on the deformation energy with nodal coordinates as variables of the optimization instead of the lengths or dimensions of the elements of the mechanism, as it was usual in this method until now because the parameters (lengths) were linear in that way. However, in that work no geometric constraint was implemented in the algorithm yet.

In the current manuscript, though, the implementation of the geometric constraints is the innovation of the paper, and the nodal coordinates are used as the optimization variables. This yields the advantage whereby the assembly configuration of the mechanism is considered in the optimization process, and is not fixed or predefined by the initial position of the links. The dimensions are defined as functions of the initial nodal coordinates, because of which neither the lengths nor the constraints are no longer linear. This approach leads to the basic drawback whereby the derivatives are more complex. Nevertheless, once they have been calculated, the solving algorithm searches in a wider workspace for the final mechanism. The algorithm is then free to search for a solution among a greater number of possible minima while satisfying the constraints applied. When the variables of the optimization are the dimensions of the links, the search is confined to a local minimum of the error function.

Some of the constraints introduced define the distance between some pairs of nodes of the mechanism, and the remainder describe the length of one or more links of the mechanism [17, 63]. In this last reference, a presentation of the method followed in the actual work is made, but some novelties are added here. To summarize: (i) A comparison between both approaches (Euclidean and alternative) to define the constraints in the Lagrange multipliers method is added. (ii) A more completely developed explanation of the implementation is given in this work. (iii) Restrictions referring to links of invariable length are introduced, which are implemented directly into the synthesis problem loop, and not in the analysis loop as happens with distance between nodes. The formulation of both constraints results in very similar mathematical expressions but this is just a casualty, because both formulations are also similar, being the difference the type of variables optimized. (iv) The application, implementation, and performance evaluation of an algorithm to impose these two geometric constraints is the major contribution of this paper.

In this paper the kinematic analysis is carried out using analytical second order derivatives, leading thus to a second order method. In the other hand, the synthesis problem, which is run fewer times, only uses approximate derivatives by central differences, leading thus to a method between linear and quadratic order.

This method, in the form presented in this paper, is limited to kinematic synthesis of mechanisms. It is also reduced to a set of precision points defined by the user and, thus, does not consider order, branch or circuit problems, although one can define the set of synthesis points to somehow coerce the method to consider them. When looking for the synthesis solution, the algorithm presented goes from one prescribed point to the next, so that it takes as initial position for solving one precision position problem, the solution of the previous position. Therefore, it is more probable that the process maintain the same assembly of the links, which is desirable, if the positions are near enough. So the branch, order and circuit are stablished by following the same assembly as in the previous prescribed point, not allowing branch or circuit changes. Besides, as the method is based on a SQP approach, the use of this along with the neighborhood of the precision points is what reduces the chance of configuration assembly changes. This leads to an algorithm less prone to assembly configuration changes, which is good given that the aim of the algorithm is to reach usable mechanisms, see [55] to [62]. In those cases where no initial guess is provided, other approaches such as heuristics can be used to generate the initial guess as, for example, genetic algorithms (see [61]).

In the recent literature, not this innovation nor similar studies have been found in any other work. The authors have searched among the last two years works related to generic methods for kinematic synthesis with geometric constraints. Unfortunately, they only found methods for dimensional synthesis on hybrid parallel mechanisms with other types of constraints including installation, joint motion limitation, joint load limitation and validity of data, as reference [64]; or constraints on lengths of the type higher than or less than, as reference [65]. Methods for kinematic synthesis for four-bar linkages were also found, as references [66–68] without constraints, and [69] who used constraints on lengths of the type higher than or less than. Another method for kinematic synthesis with path generation for a six-bar linkage Stephenson III was referred to in reference [70] without any type of constraints. Finally, there are general purpose synthesis methods for planar linkages as reference [71] where no constraint applies; or reference [72], due to the authors, which is the only one using geometric constraints of fixed length.

### 1.4 Organization of the paper

The remainder of this paper is organized as follows: In the second section, the synthesis of mechanisms by the deformed position procedure is explained. The use of the geometrical constraints of fixed distance are then presented in section three, followed by those of fixed length in section four. In sections five and six, some validation examples are provided to demonstrate the convergence of the method in several formulations. Finally, the conclusions of this study are discussed in section seven.

## 2 Synthesis of mechanisms by deformed position procedure

Before the presentation of the innovation represented by the geometric constraints implemented to the *Deformed Position Method*, the authors would like to mention a brief introduction to understand better its basis.

One of the first studies to introduce the *deformed position problem* was in the early eighties with [55], then [58] and [10] when, while studying a solution to the initial position problem, Avilés defined the following problems based on the deformation energy:

*Initial position problem* This type of problems is solved on the basis of the length of every link, the position of the ground nodes, and the input element positions. The objective function to minimize is described by the deformation energy of the mechanism, which can be written as a function of the lengths of the $b$ links:

$$\varphi = \sum_{i=1}^{b} (l_i - L_i)^2 \tag{1}$$

Or as a function of the coordinate vector $\{x\}$:

$$\varphi(\{x\}) = \sum_{i=1}^{b} \left[ (\{x\}^T [\overline{g}_i] \{x\})^{\frac{1}{2}} - L_i \right]^2 \tag{2}$$

where $l_i$ represents the length of deformed element $i$ and $L_i$ is the desired length of the same element $i$. When the length of the element $i$ is written with nodal coordinates as variables of the problem as in Eq. (2), vector $\{x\}$ collects all the coordinates of every link of the mechanism, and $[\overline{g}_i]$ is the geometric matrix of the element $i$ expanded to the dimension or size of the complete system to be coherent with the size of vector $\{x\}$. In the general case, where lengths of quite different scales appear, the use of weighting factors in Eqs. (1) and (2) might be needed. In this problem a prescribed variation is applied to the input parameter $\theta$

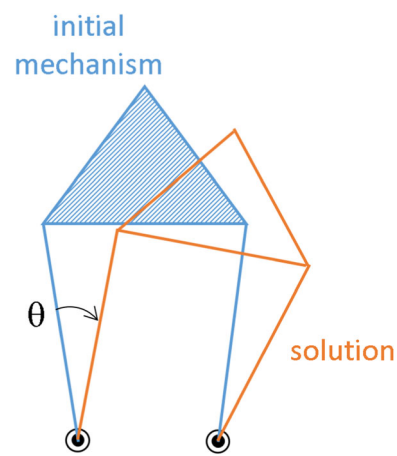

**Fig. 1** Initial position problem

 

to move the initial mechanism to the solution position as shown in Fig. 1.

*Finite displacement problem*: This problem is also called Successive Positions. In this type of problems, one seeks to determine the successive positions of the elements for a given sequence of the variation of the input parameter $\theta_i$. Thus, the displacement problem is solved as if it were the initial position problem, but taking as the initial position of each step the position solution of the previous step. See Fig. 2.

*Deformed position problem*: The function $\varphi$ defined in the initial position problem by Eqs. (1) or (2) is used here too, with the difference that the minimum value of the function is not zero but an unknown positive value near zero. It involves obtaining the position of minimum energy of the elements of a mechanism when one or some of its nodes or elements are obliged to satisfy some geometrical conditions out of the field of possible motions of a rigid solid of the mechanism. The mechanism is considered composed of deformable elements exhibiting linear elastic behavior. Thus, the error function used to measure the fitness of the mechanism with respect to a stated problem is the sum of the deformation energy stored in every link of the mechanism, which will be in the best case near zero, but that depends on the ease or possibility to fulfil the required conditions. See Fig. 3. Here, it is desired that the coupler point P reaches the prescribed position Q with the least deformation energy of the whole mechanism.

*Static equilibrium problem*: It is defined to obtain the equilibrium position of a mechanism containing elastic elements when undergoing certain external actions.

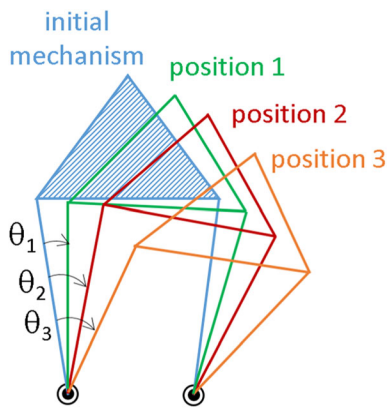

**Fig. 2** Finite displacement problem

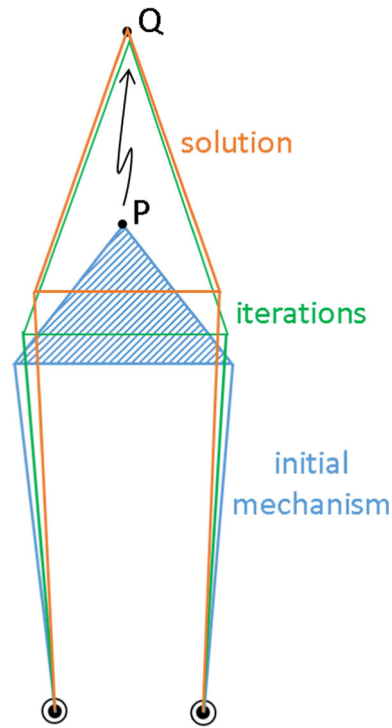

**Fig. 3** Deformed position problem

With respect to the *optimal synthesis*, this is proposed as a function of the deformed position problem. This problem is solved for each $j$ prescribed precision point of the synthesis, and the set of these points describe the desired path for the mechanism to follow. Weighting parameters $w_j$ to give more importance to some positions can be used:

$$\Phi(\{L\}) = \sum_{j=1}^{s}\left[w_j\sum_{i=1}^{b}(l_{ij} - L_i)^2\right] \qquad (3)$$

For each one of the $j$ synthesis positions, $l_{ij}$ represents the length of the deformed element $i$. Observe that this equation can also be represented dependent on the coordinates vector as it was done in Eq. (2).

All these problems presented in this section are set out as dependent on the deformed position problem. In this work the weighting parameters for every prescribed precision point have been set to the unit so that the method treats all prescribed precision positions with equal importance.

## 3 Using geometrical constraint of fixed distance between nodes as a requirement in the deformed position function

As exposed before, two new features are presented in this paper. The first one allows to introduce a new type of requirement in a synthesis step, which consist on a required distance among two nodes. This can be useful, for example, if one wants a gripper to keep a determinate distance to hold an item. It is also useful if one wants a given spring to exert a determinate force. In Fig. 4 an example is provided. Here a prescribed distance of 4.8 units is sought between nodes E and F, and in the moment of the picture one can see that the distance is a bit larger. It can also be possible to apply more than one distance restriction to the mechanism in the kinematic analysis. As the kinematic analysis of deformed position problem is solved in every synthesis step, one can expect a constraint in each step, and this constraint can be always the same between the two involved nodes or can be different for each step. This type of constraint would be linear if one uses length variables in the analysis problem; but in this work, initial nodal coordinates are used. This translates into nonlinear constraints that can be expressed as follows, where $r$ is the total number of constraints of this type in the problem:

$$c_z(\{x\}) \equiv d_z(\{x\}) - D_z = 0 \text{ for } z = 1, 2, \ldots, r \quad (4)$$

For each one of the $z$ constraints, in this expression $d_z(\{x\})$ is the real distance between the constrained nodes, say node $A$ and node $B$, for the current position, whereas $D_z$ is the restricted value for the fixed distance. The position of the mechanism changes through the necessary iterations of the kinematic analysis of the deformed position problem, to approximate $d_z(\{x\})$ to $D_z$ as near as possible. This value of $D_z$ is fixed for each prescribed precision point, i.e., if the synthesis is defined by $S$ precision points, one can define just one or $S$ different values of $D_z$. Therefore, the optimization problem could be expressed by the energy function with point-to-point distance restrictions as:

$$\min_x f(\{x\}) = \sum_{j=1}^s \sum_{i=1}^b \left[ l_{ij}(\{x\}) - L_{ij}(\{X\}) \right]^2 \quad (5)$$
$$\text{subject to: } c_z(\{x\}) = 0$$

If one writes the following:

$$\tilde{f}(\{x\}) = \sum_{i=1}^b \left[ l_{ij}(\{x\}) - L_{ij}(\{X\}_j) \right]^2 \quad (6)$$

Now, this is the function for the kinematic analysis for a synthesis position $j$. The floating nodal coordinates $\{x\}$ are the only variables, as for each synthesis position $j$ the ground nodes $\{X\}_j$ are immovable and are the result from the previous synthesis iteration. Consequently, the optimization problem can also be written as:

$$\min_x f(\{x\}) = \sum_{j=1}^s \tilde{f}(\{x\}_j) \quad (7)$$
$$\text{subject to: } c_z(\{x\}) = 0$$

For each synthesis position $j$, a kinematic position problem is carried out. Therefore, every synthesis step $j$ a set of nodal coordinates is achieved, $\{x\}_j$.

To introduce this type of constraint to the problem, as explained right before Sect. 1.2, penalty functions, the augmented Lagrangian, or Lagrange multipliers can be used.

If one uses Lagrange multipliers, the problem is defined by the next expression of the augmented error function:

$$L(\{x\}, \{\lambda\}) = \tilde{f}(\{x\}) - \sum_{z=1}^r \lambda_z c_z(\{x\}) \quad (8)$$

In this expression, $\tilde{f}(\{x\})$ is the error function of the deformed position problem in Eq. (6), and the second term corresponds to the sum of every constraint $c_z$ in Eq. (4), multiplicated by the Lagrange multipliers $\lambda_z$.

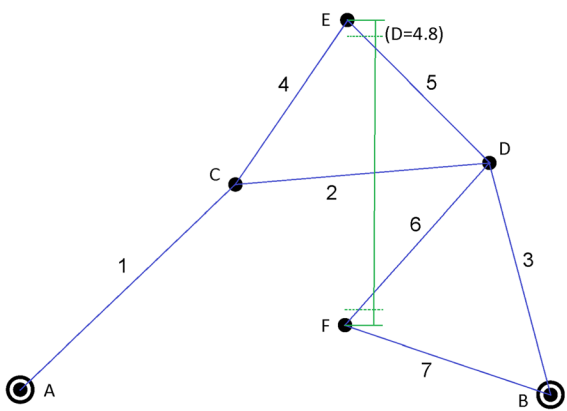

**Fig. 4** Mechanism with fixed distance of $D_z = 4.8$ units between two floating nodes

Moreover, two formulations are proposed to define the constraints: If $A(x_a, y_a)$ and $B(x_b, y_b)$ are the constrained nodes of the mechanism, then:

1. Euclidean distances:

$$c_z(\{x\}) = d_z(\{x\}) - D_z =$$
$$= \sqrt{(x_{b_z} - x_{a_z})^2 + (y_{b_z} - y_{a_z})^2} - D_z = 0 \tag{9}$$

2. Alternative formulation (difference of squared distances):

$$\tilde{c}_z(\{x\}) = d_z(\{x\})^2 - D_z^2 =$$
$$= (x_{b_z} - x_{a_z})^2 + (y_{b_z} - y_{a_z})^2 - D_z^2 = 0 \tag{10}$$

where $D_z$ in both cases is the constrained distance $z$.

### 3.1 Analytical development of the Euclidean distance formulation

To apply the above constraints, first and second derivatives with respect to $\{x\}$ are needed. This development is valid for the three aforementioned methods.

(a) Derivative of constraint $c_z(\{x\})$, Eq. (9), with respect to the variables $\{x\}$, gradient vector:

$$\left\{\frac{\partial c_z(\{x\})}{\partial \{x\}}\right\} = \{g_z(\{x\})\} = \frac{1}{d_z(\{x\})} \begin{Bmatrix} -(x_{b_z} - x_{a_z}) \\ -(y_{b_z} - y_{a_z}) \\ (x_{b_z} - x_{a_z}) \\ (y_{b_z} - y_{a_z}) \end{Bmatrix} \tag{11}$$

(b) Derivative of constraint $c_z(\{x\})$ with respect to the variables $\{x\}$ twice, Hessian matrix:

$$\left[\frac{\partial^2 c_z(\{x\})}{\partial \{x\}^2}\right] = [h_z(\{x\})]$$
$$= \frac{1}{d_z(\{x\})}$$
$$\begin{bmatrix} 1 - l_{xz}^2 & -l_{xz} \cdot l_{yz} & -1 + l_{xz}^2 & l_{xz} \cdot l_{yz} \\ -l_{xz} \cdot l_{yz} & 1 - l_{yz}^2 & l_{xz} \cdot l_{yz} & -1 + l_{yz}^2 \\ -1 + l_{xz}^2 & l_{xz} \cdot l_{yz} & 1 - l_{xz}^2 & -l_{xz} \cdot l_{yz} \\ l_{xz} \cdot l_{yz} & -1 + l_{yz}^2 & -l_{xz} \cdot l_{yz} & 1 - l_{yz}^2 \end{bmatrix} \tag{12}$$

where $l_{xz} = \dfrac{x_{b_z} - x_{a_z}}{d_z(\{x\})}$ and $l_{yz} = \dfrac{y_{b_z} - y_{a_z}}{d_z(\{x\})}$

The derivatives are more complex than in the alternative formulation, and are terms dependent on distance $d_z(\{x\})$ as it can be seen on Eqs. (11) and (12). It is important to highlight the case where this distance becomes null, either instantly in an iteration or because it needs to be zero. In this situation, the derivatives results in an indeterminate outcome, 0/0. Under these circumstances, the gradient can be forced to be that of the alternative formulation (Eq. (13)) and the Hessian can be forced to be the unit matrix, or the same resulting matrix for the alternative formulation (Eq. (14)).

### 3.2 Analytical development of the alternative formulation

Once again, this development is valid for the three methods mentioned in Sect. 1.1 for nonlinear equality constraints.

(a) Derivative of constraint $\tilde{c}_z(\{x\})$, Eq. (10), with respect to variables $\{x\}$, gradient vector:

$$\left\{\frac{\partial \tilde{c}_z(\{x\})}{\partial \{x\}}\right\} = \{\tilde{g}_z(\{x\})\} = 2 \begin{Bmatrix} -(x_{b_z} - x_{a_z}) \\ -(y_{b_z} - y_{a_z}) \\ (x_{b_z} - x_{a_z}) \\ (y_{b_z} - y_{a_z}) \end{Bmatrix} \tag{13}$$

(b) Derivative of constraint $\tilde{c}_z(\{x\})$ with respect to variables $\{x\}$ twice, Hessian matrix:

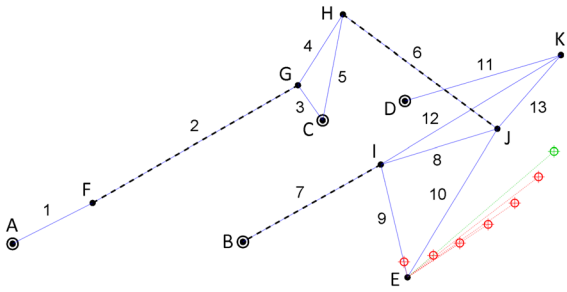

**Fig. 5** Mechanism with fixed lengths in three links

$$\left[\frac{\partial^2 \tilde{c}_z(\{x\})}{\partial \{x\}^2}\right] = [\tilde{h}_z(\{x\})] = \begin{bmatrix} 2 & 0 & -2 & 0 \\ 0 & 2 & 0 & -2 \\ -2 & 0 & 2 & 0 \\ 0 & -2 & 0 & 2 \end{bmatrix}$$

$$(14)$$

## 4 Using geometrical constraints of fixed length in the synthesis function

Another type of constraint used in this work is the fixing of the length of certain elements of the mechanism, which are introduced in the synthesis problem. Figure 5 shows an example. This implies inclusion of nonlinear constraints in the synthesis problem. These constraints are linear when expressed in dimensions but nonlinear when expressed under initial nodal coordinates. They can be written as follows, where $r$ is the total number of constraints of the type included in the problem:

$$c_t(\{X\}) \equiv L_t(\{X\}) - \bar{L}_t = 0 \text{ for } t = 1, 2, \ldots, r$$

In the above expression, $L_t(\{X\})$ represents the real length of the constrained element $t$, say between ends $A$ and $B$, and $\bar{L}_t$ is the fixed length for the element $t$. As this constraint is included in the synthesis problem, it changes for each step of the iteration as the ground nodal coordinates $\{X\}$ change and consequently, $L_t(\{X\})$ changes, and approaches the value desired, $\bar{L}_t$. This value of $\bar{L}_t$ is fixed for the overall problem, that is, it is the length desired for the element, and can only be one during the problem. Thus, the optimization problem could be expressed by the energy function with fixed length resctrictions as in Eq. (15):

$$\min_X f(X) = \sum_{j=1}^{s} \sum_{i=1}^{b} \left[l_{ij}(x) - L_{ij}(X)\right]^2$$

$$\text{subject to: } c_t(X) = 0 \quad (15)$$

Once the kinematic analysis has been carried out in the synthesis step $j$, $\{x\}$ now are the nodal coordinates of every floating node acquired in the synthesis position $j$, that is $\{x\}_j$, and will not change during the synthesis problem; and $\{X\}$ belong to the ground nodes coordinates in the same synthesis position $j$, which are the optimization variables for this position $j$, then one can call them $\{X\}_j$.

Accordingly, one can also write for the optimization problem the following expression of Eq. (16):

$$\min_X f(\{X\}) = \sum_{j=1}^{s} \sum_{i=1}^{b} \left[l_{ij}(\{x\}_j) - L_{ij}(\{X\}_j)\right]^2$$

$$\text{subject to: } c_t(\{X\}) = 0 \quad (16)$$

One example for this type of constraint is given in Fig. 5. Here a packaging mechanism is shown and for functional reasons some of its elements have been restricted. These constrained lengths are marked by discontinuous black lines: bar 2 between nodes F and G, bar 6 between nodes H and J, and bar 7 between nodes B and I. These lengths are maintained constant along the synthesis, while node E searches the 8 prescribed positions as near as possible marked by the circles. The position represented in green means the beginning of the desired path.

For this type of constraint, different approaches can be used, such as penalty functions, Lagrange multipliers, and Augmented Lagrangian. One can get an equivalent equation to that from (8):

$$L(\{X\}, \{\lambda\}) = f(\{X\}) - \sum_{z=1}^{r} \lambda_z c_z(\{X\}) \quad (17)$$

Where $f(\{X\})$ is now the error function of the synthesis problem in Eq. (15).

Two formulations can also be proposed, as in previous section: if $A(X_a, Y_a)$ and $B(X_b, Y_b)$ are the end nodes of the constrained link,

1. Euclidean length:

$$c_t(\{X\}) = L_t(\{X\}) - \tilde{L}_t =$$
$$= \sqrt{(X_{b_t} - X_{a_t})^2 + (Y_{b_t} - Y_{a_t})^2} - \tilde{L}_t = 0 \tag{18}$$

2. Alternative formulation (difference of squared distances):

$$\tilde{c}_t(\{X\}) = L_t(\{X\})^2 - \tilde{L}_t^2 =$$
$$= (X_{b_t} - X_{a_t})^2 + (Y_{b_t} - Y_{a_t})^2 - \tilde{L}_t^2 = 0 \tag{19}$$

where $\tilde{L}_t$ is the constrained length $t$.

Now, the development of the two approaches is very similar to that of the fixed distance described in the previous Sect. 3.1, with the only difference that before when $d_z(\{x\})$, $x_{b_z}$ and $x_{a_z}$ were used, now $L_t(\{X\})$, $X_{b_t}$ and $X_{a_t}$ will be used. Thus, only the Euclidean length will be advanced this time.

### 4.1 Analytical development of the Euclidean length formulation

To apply these constraints to the synthesis problem, first and second derivatives with respect to $\{X\}$ are needed. This development is valid for the aforementioned three methods.

(a) Derivative of constraint $c_t$ with respect to variables $\{X\}$, gradient vector:

$$\left\{\frac{\partial c_t(\{X\})}{\partial \{X\}}\right\} = \{g_t(\{X\})\}$$
$$= \frac{1}{L_t(\{X\})} \left\{ \begin{array}{c} -(X_{b_t} - X_{a_t}) \\ -(Y_{b_t} - Y_{a_t}) \\ (X_{b_t} - X_{a_t}) \\ (Y_{b_t} - Y_{a_t}) \end{array} \right\} \tag{20}$$

(b) Derivative of constraint $c_t$ with respect to variables $\{X\}$ twice, Hessian matrix:

$$\left[\frac{\partial^2 c_t(\{X\})}{\partial \{X\}^2}\right] = [h_t(\{X\})]$$
$$= \frac{1}{L_t(\{X\})}$$
$$\begin{bmatrix} 1 - l_{Xt}^2 & -l_{Xt} \cdot l_{Yt} & -1 + l_{Xt}^2 & l_{Xt} \cdot l_{Yt} \\ -l_{Xt} \cdot l_{Yt} & 1 - l_{Yt}^2 & l_{Xt} \cdot l_{Yt} & -1 + l_{Yt}^2 \\ -1 + l_{Xt}^2 & l_{Xt} \cdot l_{Yt} & 1 - l_{Xt}^2 & -l_{Xt} \cdot l_{Yt} \\ l_{Xt} \cdot l_{Yt} & -1 + l_{Yt}^2 & -l_{Xt} \cdot l_{Yt} & 1 - l_{Yt}^2 \end{bmatrix} \tag{21}$$

where $l_{Xt} = \dfrac{X_{b_t} - X_{a_t}}{L_t(\{X\})}$ and $l_{Yt} = \dfrac{Y_{b_t} - Y_{a_t}}{L_t(\{X\})}$

Note that the development of the derivatives is very similar to that of the fixed distance, with the only difference being that in this case, the design variables are the initial coordinates $\{X\}$ instead of nodal coordinates $\{x\}$ of the deformed problem.

When working with this optimization problem without enough constraints, it is well known that it might be convex or semidefinite ([62]). The Hessian matrix in case of the convex problem is positive definite, and the algorithm converges to a local minimum. However, when introducing nonlinear constraints by means of Lagrange multipliers or the augmented Lagrangian, the problem is no longer convex, as it always converges to a saddle point.

The algorithm used for synthesis optimization is based on the SQP method such that it solves sequentially a series of quadratic problems, using for each one of them the position of the solution of the previous problem. This method leads to a system known as Karush–Kuhn–Tucker [73–75].

$$\begin{bmatrix} [H_L]_{x^2} & [H_L]_{x\lambda} \\ [H_L]_{x\lambda} & [H_L]_{\lambda^2} \end{bmatrix} \left\{ \begin{array}{c} \Delta\{x\} \\ \Delta\{\lambda\} \end{array} \right\} = -\left\{ \begin{array}{c} \{G_L\}_x \\ \{G_L\}_\lambda \end{array} \right\} \tag{22}$$

One can observe that this expression is of the type $[H]\{\Delta x\} = \{G\}$, that is, the Hessian matrix of the system multiplied by the vector of variables of the system equal to minus the gradient of the system.

The first element in the matrix $[H_L]_{x^2}$ is the total Hessian of the augmented error function, when derivating with respect to the nodal coordinates twice. To simplify the expression it can be called $[H]$. The second element of the matrix $[H_L]_{x\lambda}$ is formed by the derivatives of the constraints with respect to $\{x\}$. Take into account that the derivatives of the augmented error function with respect to the Lagrange multipliers $\{\lambda\}$ give as a result the constraints functions. To simplify the expression it can be called $[E]$. This matrix is not square, therefore, the matrix and its transpose will be needed to fill the main matrix. The last element in the matrix $[H_L]_{\lambda^2}$ is the null matrix $[0]$ because when deriving twice with respect to the multipliers, as there is no term with $\lambda$ squared, it will be zero. In the vector of the variables, one can find first the nodal coordinates that have to be found to achieve the solution, and also the Lagrange multipliers themselves, because they have to be found to define the system. Finally, to the right of the expression, there are two elements in the vector. The first element $\{G_L\}_x$ consists of the total gradient of the augmented error function, when derivating with respect to the nodal coordinates. To simplify the expression this will be called $\{G\}$. Lastly, the second element $\{G_L\}_\lambda$ is composed by all the constraints. In case of use of the Lagrange multipliers, the gradient with respect to the $\{\lambda\}$ will be properly a vector conformed by the constraints.

And considering the above, the expression (22) can also be written as (23):

$$\begin{bmatrix} [H] & [E]^T \\ [E] & [0] \end{bmatrix} \begin{Bmatrix} \Delta\{x\} \\ \Delta\{\lambda\} \end{Bmatrix} = -\begin{Bmatrix} \{G\} \\ \{R\} \end{Bmatrix} \qquad (23)$$

In this system the total dimension considers all the nodal coordinates (2 times the number of nodes of the mechanism) plus all the Lagrange multipliers necessary to introduce the constraints in the system (one multiplier per constraint). To solve the desired parameters $\Delta\{x\}$, the Null Space method is applied as introduced by F. de Bustos et al. in [63]. Furthermore, in order to solve saddle points on undesired extrema, one can use modified factorizations or negative curvature search algorithms.

About the computational cost of the method it can be discussed separately in relation to the analysis of the position problem and the synthesys. For the analysis, it lies in the order of 2 to $4 \times 10^{-4}$ seconds for simple topologies like four-bars, and it can be increased to about 7 to $8 \times 10^{-4}$ seconds for more complex topologies like mechanisms with several links as the one studied in example 3 in next section. For the synthesis, though, the working time increases as the process consists on several internal loops of analysis problems, one for each prescribed precision point of the synthesis. Therefore, the computational cost will be approximately in the order of the seconds. Obviously, this cost is directly proportional to the number of prescribed points as it involves more analysis problems to be solved, and additionally proportional too, to the number of nodes to be positioned.

## 5 Examples to verify requirements of fixed distance

In spite of the wide range of methods for synthesis existing in the bibliography, the features introduced in this paper have been seldom studied. Thus, the authors have not found examples with geometric restrictions of the type presented here in order to be able to compare results. Therefore, to prove the appropriateness of the behavior of the algorithm, ad-hoc examples have been developed. These cover different types of problems, such as:

1. complexity of the topology taking into account the number of links of the mechanism;
2. distance between the initial guess and the final solution, where it is an initial position problem if it is large and a successive positions problem if it is small; and
3. whether the final solution is deformed, i.e., if the position of the solution can be attained by deforming the mechanism.

The objectives to prove were the following:

1. Verify the convergence of the function of the deformation energy as well as the fulfillment of nonlinear constraints.
2. Verify the two formulations of the nonlinear constraints (Euclidean and alternative).
3. Verify the behavior of the function based on the initial nodal coordinates.
4. Verify modifications to the optimization system.

For the examples shown, a random initial guess is tried and the method looks for a solution. Therefore, the chosen position will influence on the solution found. If the initial guess of the algorithm is far from the final solution of the position problem, the way to the local minima for the Euclidean and Alternative approaches is slightly different and this can lead to some other local minima due to a distinct configuration of the final multicore capabilities have not been used. This has been done due to the fact that the cost of creating the thread is orders of magnitude higher than the cost of solving the linear systems.

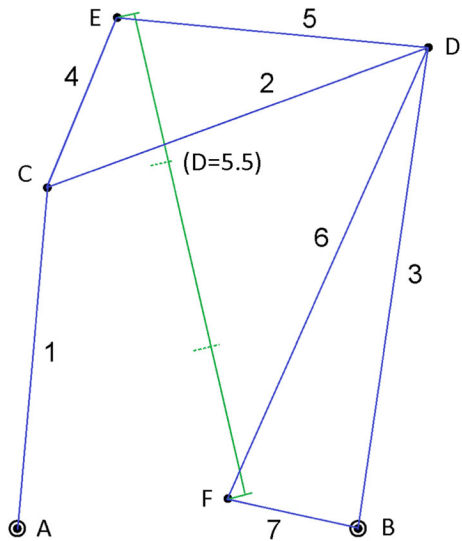

**Fig. 6** Initial mechanism of four-link topology, example 1

**Table 1** Nodal coordinates of example 1 in initial and final configurations

|         | $x_A$  | $y_A$  | $x_B$ | $y_B$ | $x_C$ | $y_C$ |
|---------|--------|--------|-------|-------|-------|-------|
| Initial | 0.00   | 0.00   | 9.84  | 0.00  | 0.87  | 9.84  |
| Final   | 0.00   | 0.00   | 9.84  | 0.00  | 9.60  | 2.33  |

|         | $x_D$ | $y_D$ | $x_E$ | $y_E$ | $x_F$ | $y_F$ |
|---------|-------|-------|-------|-------|-------|-------|
| Initial | 11.87 | 13.89 | 2.89  | 14.76 | 6.08  | 0.87  |
| Final   | 6.66  | 13.67 | 4.91  | 4.82  | 6.02  | −0.56 |

### 5.1 Four-link topology, initial position problem, and low deformation energy

The mechanism of example 1 shows a four-bar with restricted distance between nodes E and F in D = 5.5 units. In the initial guess (Fig. 6) this distance was far from the solution (Fig. 7), and thus needed several steps to get to the solution. Although the initial guess was far from the desired distance, the topological assembly did not change, so the needed deformation energy was not high.

Table 1 shows the coordinates of every node in the four-bar in the initial guess and in the final solution. From this table one can obtain the initial and final distance between nodes E and F to check the fulfilment of the distance restriction D.

The optimization problem here is a kinematic analysis of deformed position problem defined in Eq. (24):

$$\min_x f(\{x\}) = \sum_{i=1}^{7} [l_i(\{x\}) - L_i(\{X\})]^2$$
$$\text{subject to: } c(\{x\}) = d(\{x\}) - 5.5 = 0 \tag{24}$$

Table 2 exhibits the values of the deformation energy for both approaches in the 19 or 20 iterations. Moreover, Table 3 presents the values of the error calculated for both approaches in the 16 or 21 iterations. All these results from Tables 2 and 3 have been plotted in Fig. 8 showing graphically that in solution. In all the examples presented, both approaches lead to numerically equal results. For this reason, only one solution for each example is displayed. Some representative cases of a kinematic analysis of the deformed position problem are provided below. The computer used was an INTEL XEON E5-2637 v3 @3.5GHz. As the matrices used are of small size, a single thread has been used so example 1, both the convergence of energy and the error was faster in the Euclidean approach, although both approaches tended to oscillate, more in the constraint error than in the deformation energy. One can observe that this oscillation is not always a negative aspect, as it drives to a final convergence at a considerably low iteration amount as 21. Total time was of 0.0004 s for the fastest approach.

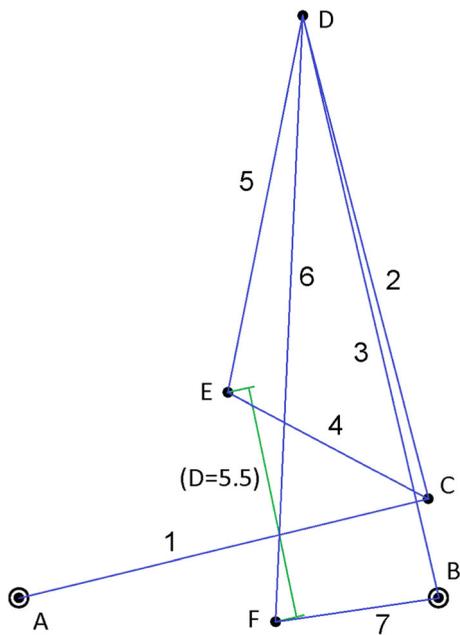

**Fig. 7** Solution mechanism of four-link topology, example 1

**Table 2** Resulting values of energy during the iterations for example 1 in both Euclidean and alternative approaches

| Iter | Euclidean fitness | Log fitness | Alternative fitness | Log fitness |
|------|-------------------|-------------|---------------------|-------------|
| 0 | 0.00 | | 0.00 | |
| 1 | 166.54 | 2.22 | 48.32 | 1.68 |
| 2 | 1.38 | 0.14 | 0.44 | − 0.36 |
| 3 | 0.94 | − 0.03 | 0.85 | − 0.07 |
| 4 | 0.50 | − 0.30 | 0.49 | − 0.31 |
| 5 | 0.71 | − 0.15 | 0.35 | − 0.45 |
| 6 | 0.23 | − 0.64 | 0.26 | − 0.58 |
| 7 | 2.63 | 0.42 | 0.20 | − 0.69 |
| 8 | 0.06 | − 1.25 | 0.14 | − 0.85 |
| 9 | 1.27 | 0.11 | 0.29 | − 0.54 |
| 10 | 0.01 | − 2.25 | 0.06 | − 1.25 |
| 11 | 0.06 | − 1.24 | 0.98 | − 0.01 |
| 12 | 1.08E−04 | − 3.97 | 0.01 | − 2.15 |
| 13 | 4.24E−05 | − 4.37 | 0.08 | − 1.07 |
| 14 | 1.44E−08 | − 7.84 | 1.69E−04 | − 3.77 |
| 15 | 8.71E−13 | − 12.06 | 1.04E−04 | − 3.98 |
| 16 | 4.58E−23 | − 22.34 | 5.38E−08 | − 7.27 |
| 17 | 6.51E−30 | − 29.19 | 1.20E−11 | − 10.92 |
| 18 | 1.26E−29 | − 28.89 | 7.47E−21 | − 20.13 |
| 19 | 1.97E−29 | − 28.71 | 2.62E−29 | − 28.58 |
| 20 | | | 1.97E−29 | − 28.71 |

### 5.2 Four-link topology, successive position problem, and low deformation energy

In this successive problem, example 2 shows a four-bar mechanism which corresponds to the solution of example 1, but where nodes have been moved a small quantity from the solution position. Now the initial guess (Fig. 9) is quite near the solution position (Fig. 10), and this is why it is a succesive position problem. Distance between nodes E and F is restricted to D = 5.5 units as before.

Table 4 manifests the nodal coordinates of the initial guess and the final solution, and hereby one can calculate the initial and final distance between restricted nodes to verify the fulfilment of the condition.

The optimization problem here is a kinematic analysis of deformed position problem defined in Eq. (25), which is of course the same as in example 1:

$$\min_x f(\{x\}) = \sum_{i=1}^{7} [l_i(\{x\}) - L_i(\{X\})]^2 \qquad (25)$$

$$\text{subject to: } c(\{x\}) = d(\{x\}) - 5.5 = 0$$

Table 5 displays the values of the deformation energy in every 10 iterations, for both approaches. Table 6, in its turn, manifests the values of the error in the constraint in every 9 iterations, for both approaches.

The results in Tables 5 and 6 are plotted in Fig. 11. It can be observed that the convergence of both the energy and the error was faster than in example 1, as expected because the initial guess was nearer the final position. The behavior of both approaches shows a high similarity between them to the point that the lines in the plot are confused, meaning that for succesive positions problems, the alternative approach is as good as the Euclidean one. The example took 0.0002 s to be

**Table 3** Resulting values of error during the iterations for example 1 in both Euclidean and alternative approaches

| Iter | Euclidean error | Log error | Alternative error | Log error |
| --- | --- | --- | --- | --- |
| 0 | 8.75 | 0.94 | 172.86 | 2.24 |
| 1 | 8.47 | 0.93 | 115.95 | 2.06 |
| 2 | 0.24 | − 0.62 | 24.95 | 1.40 |
| 3 | 2.55E−04 | − 3.59 | 5.26 | 0.72 |
| 4 | 4.40E−04 | − 3.36 | 0.30 | − 0.53 |
| 5 | 0.13 | − 0.87 | 0.43 | − 0.37 |
| 6 | 7.51E−05 | − 4.12 | 0.13 | − 0.89 |
| 7 | 0.25 | − 0.59 | 0.41 | − 0.39 |
| 8 | 6.40E−04 | − 3.19 | 0.04 | − 1.45 |
| 9 | 0.14 | − 0.85 | 0.79 | − 0.10 |
| 10 | 2.58E−04 | − 3.59 | 5.46E−03 | − 2.26 |
| 11 | 0.03 | − 1.59 | 1.39 | 0.14 |
| 12 | 8.59E−06 | − 5.07 | 0.02 | − 1.76 |
| 13 | 6.40E−04 | − 3.19 | 0.35 | − 0.46 |
| 14 | 6.22E−08 | − 7.21 | 1.16E−03 | − 2.94 |
| 15 | 8.99E−08 | − 7.05 | 1.10E−02 | − 1.96 |
| 16 | 1.51E−14 | − 13.82 | 1.82E−06 | − 5.74 |
| 17 | | | 3.68E−06 | − 5.43 |
| 18 | | | 2.07E−12 | − 11.68 |
| 19 | | | 7.11E−15 | − 14.15 |
| 20 | | | 3.55E−15 | − 14.45 |
| 21 | | | 7.11E−15 | − 14.15 |

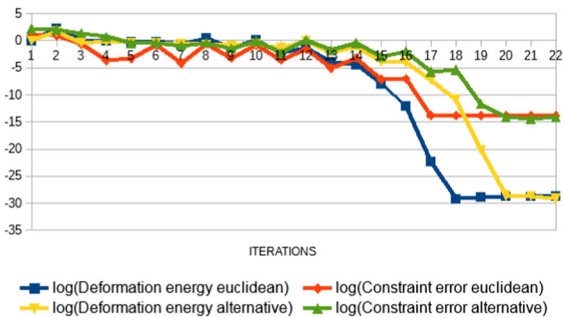

**Fig. 8** Evolution of the deformation energy, and of the error in the constraint with both approaches, example 1. The values have been displayed in logarithmic scale

solved for both approaches with negligible error less than 1%.

### 5.3 Topology with more links, initial position problem, and low deformation energy

In example 3 a mechanism with 17 elements is presented (see Fig. 12). Here distance between nodes

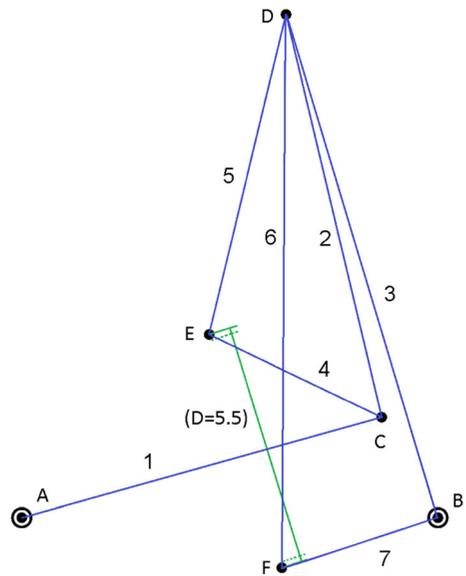

**Fig. 9** Initial mechanism of four-link topology, example 2

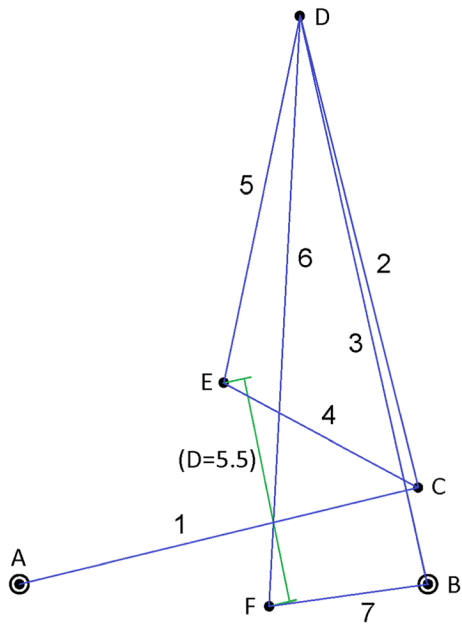

**Fig. 10** Solution mechanism of four-link topology, example 2

K and L was fixed to D=5 units. The resulting mechanism was that shown in Fig. 13.

The optimization problem was depicted by Eq. (26) for the corresponding kinematic analysis:

$$\min_x f(\{x\}) = \sum_{i=1}^{17}[l_i(\{x\}) - L_i(\{X\})]^2$$

$$\text{subject to: } c(\{x\}) = d(\{x\}) - 5.0 = 0$$

(26)

Table 7 indicates the nodal coordinates of the mechanism in example 3 in the initial and final configurations. If one centrates in nodes K and L, it is possible to calculate distance between those mentioned nodes and verify that in the solution the constraint is fulfilled.

These results of example 3 are shown in Fig. 14. Observe that they were very similar to those found in the four-links topology, with the particularity that this case incurred higher computational cost as it had a larger number of nodes to position. Considering the complexity of the topology, the convergence speed is noteworthy. Additionaly, one can remark as in the plot of example 2, that the behavior of the pair of approaches is so alike that the evolution lines on the plot get unclear.

Results of deformation energy and of error in the distance constraint are exhibited in Tables 8 and 9 respectively. There one can see that convergence was achived in 19 iterations for both approaches.

As a conclusion of this study, on the one hand, the clearer physical interpretation of the Euclidean approach should be noted (difference in distances); but on the other hand, the alternative approach

**Table 4** Nodal coordinates of example 2 in initial and final configurations

|         | $x_A$ | $y_A$ | $x_B$ | $y_B$ | $x_C$ | $y_C$ |
|---------|-------|-------|-------|-------|-------|-------|
| Initial | 0.00  | 0.00  | 9.84  | 0.00  | 8.51  | 2.38  |
| Final   | 0.00  | 0.00  | 9.84  | 0.00  | 9.60  | 2.33  |

|         | $x_D$ | $y_D$  | $x_E$ | $y_E$ | $x_F$ | $y_F$  |
|---------|-------|--------|-------|-------|-------|--------|
| Initial | 6.25  | 11.93  | 4.43  | 4.34  | 6.15  | −1.19  |
| Final   | 6.75  | 13.69  | 4.92  | 4.86  | 6.02  | −0.53  |

**Table 5** Resulting values of energy during the iterations for example 2 in both Euclidean and alternative approaches

| Iter | Euclidean fitness | Log fitness | Alternative fitness | Log fitness |
|------|-------------------|-------------|---------------------|-------------|
| 0    | 10.56             | 1.02        | 10.56               | 1.02        |
| 1    | 0.10              | −0.98       | 0.10                | −0.98       |
| 2    | 2.64E−03          | −2.58       | 2.56E−03            | −2.59       |
| 3    | 3.77E−05          | −4.42       | 3.27E−05            | −4.49       |
| 4    | 2.38E−06          | −5.62       | 1.74E−06            | −5.76       |
| 5    | 1.15E−08          | −7.94       | 6.59E−09            | −8.18       |
| 6    | 5.49E−13          | −12.26      | 1.84E−13            | −12.74      |
| 7    | 6.31E−22          | −21.20      | 7.64E−23            | −22.12      |
| 8    | 1.05E−29          | −28.98      | 7.10E−30            | −29.15      |
| 9    | 3.94E−30          | −29.40      | 3.94E−30            | −29.40      |
| 10   | 4.14E−30          | −29.38      | 1.10E−29            | −28.96      |

| Iter | Euclidean Error | Log Error | Alternative Error | Log Error |
|---|---|---|---|---|
| 0 | 0.29 | −0.54 | 3.29 | 0.52 |
| 1 | 0.09 | −1.07 | 1.03 | 0.01 |
| 2 | 0.03 | −1.55 | 0.31 | −0.51 |
| 3 | 9.75E−04 | −3.01 | 0.01 | −1.95 |
| 4 | 1.28E−04 | −3.89 | 1.19E−03 | −2.92 |
| 5 | 2.43E−06 | −5.61 | 2.13E−05 | −4.67 |
| 6 | 6.40E−08 | −7.19 | 4.03E−07 | −6.39 |
| 7 | 5.33E−13 | −12.27 | 2.15E−12 | −11.67 |
| 8 | 8.88E−16 | −15.05 | 1.07E−14 | −13.97 |
| 9 | 0.00 | | 1.42E−14 | −13.85 |

**Table 6** Resulting values of error during the iterations for example 2 in both Euclidean and alternative approaches

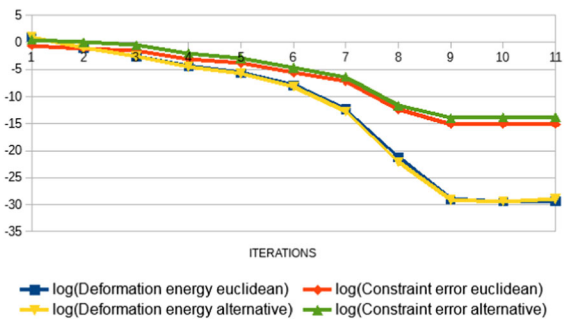

**Fig. 11** Evolution of the deformation energy, and of the error in the constraint with both approaches, example 2. The values have been displayed in logarithmic scale

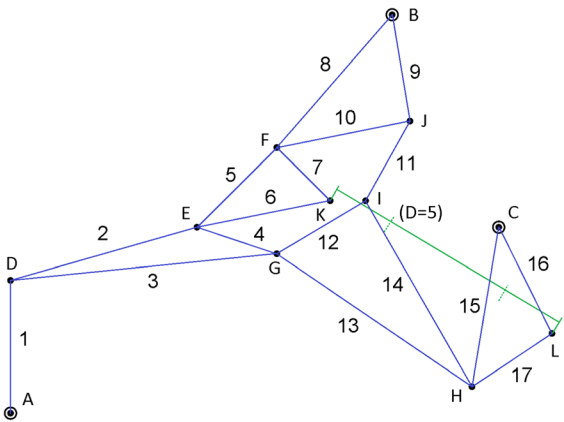

**Fig. 12** Initial mechanism of topology with more links, example 3

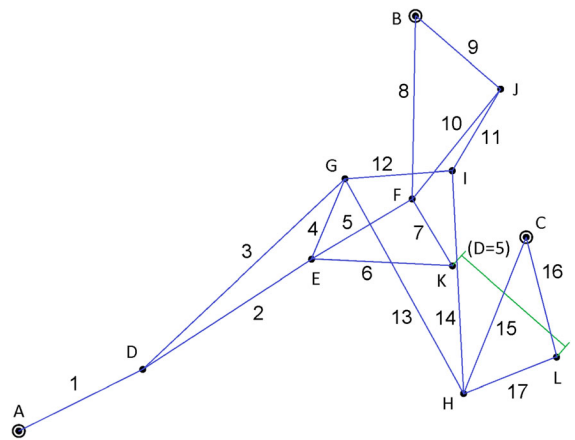

**Fig. 13** Solution mechanism of topology with more links, example 3

fashion. Therefore, either of the two approaches studied is adequate for the optimization algorithm. This example was solved in 0.00075 s in either approaches with negligible error below the 1%.

Nevertheless, as both approaches work similarly, in the next section only the Euclidean approach will be used.

## 6 Examples to confirm geometrical constraints of fixed length

Depending on the initial guess applied to the optimization algorithm, the achieved solutions can result very different. However, the method has shown to be reliable due to the fact that there is always a solution; it may happen that it be a local optimum, which not always will be a bad result as it could be interesting to

(difference of squared distances) was easier to work with, simpler when the derivatives need to be calculated, and the Hessian results in a more uncoupled matrix. With both approaches, similar results were obtained with the problem behaving in a similar find a solution similar to the initial mechanism because it may show advantages not considered in the optimization.

**Table 7** Nodal coordinates of example 3 in initial and final configurations

| | $x_A$ | $y_A$ | $x_B$ | $y_B$ | $x_C$ | $y_C$ | $x_D$ | $y_D$ |
|---|---|---|---|---|---|---|---|---|
| Initial | 0.00 | 0.00 | 14.33 | 15.00 | 18.33 | 7.00 | 0.00 | 5.00 |
| Final | 0.00 | 0.00 | 14.33 | 15.00 | 18.33 | 7.00 | 4.48 | 2.22 |

| | $x_E$ | $y_E$ | $x_F$ | $y_F$ | $x_G$ | $y_G$ | $x_H$ | $y_H$ |
|---|---|---|---|---|---|---|---|---|
| Initial | 7.00 | 7.00 | 10.00 | 10.00 | 10.00 | 6.00 | 17.33 | 1.00 |
| Final | 10.58 | 6.20 | 14.21 | 8.39 | 11.78 | 9.12 | 16.07 | 1.35 |

| | $x_I$ | $y_I$ | $x_J$ | $y_J$ | $x_K$ | $y_K$ | $x_L$ | $y_L$ |
|---|---|---|---|---|---|---|---|---|
| Initial | 13.33 | 8.00 | 15.00 | 11.00 | 12.00 | 8.00 | 20.33 | 3.00 |
| Final | 15.66 | 9.41 | 17.41 | 12.36 | 15.67 | 5.96 | 19.43 | 2.66 |

Furthermore, the increase of the number of prescribed precision points has a considerable impact, as well as in the basic algorithm presented in the references [10, 55], and [58].

Although the method does not include a specially tailored way to maintain the assembly configuration, one can coerce the method into keeping it by introducing the precision points in the neighborhood, so each consecutive precision point is the nearby of the previous one. This reduces the probability of a configuration assembly change among precision points. This is possible thanks to the use of an exploitative method (SQP), which does not lead to great changes in the initial solution if this is near enough to the local optima. Thus, the final mechanism will be better than the initial because it approximates with higher accuracy the prescribed points while fulfilling the imposed geometric constraints, which maybe the initial one does not verify.

In the following two examples, several synthesis positions are searched for. These prescribed positions are represented by circles, and the first one is painted in green to symbolize the beginning of the path. The remaining circles are painted in red.

### 6.1 Four-link topology with one constraint

In this example of a four-bar, shown in Fig. 15, it was desirable to draw a path described by 11 prescribed precision points (see Table 11) in the shape of a straight horizontal line with node E (coupler point). Theoretically a mechanism can only fulfil a straight line with 6 prescribed points at most, but the aim of the method is to approximate the path, fulfilling the positions nearly and not exactly. During the path synthesis the element placed between nodes A and C was restricted to a fixed value of D=4 units.

The optimization problem in this example 4 could be rendered by Eq. (27):

$$\min_X f(\{X\}) = \sum_{j=1}^{11} \sum_{i=1}^{5} \left[ l_{ij}(\{x\}_j) - L_{ij}(\{X\}_j) \right]^2$$
$$\text{subject to: } c_t(\{X\}) = L_t(\{X\}) - 4.0 = 0 \tag{27}$$

Table 10 presents the nodal coordinates of every node in the mechanism, both for the initial guess 15 and the solution mechanism 16. Here one can check whether the final restricted length of the element between nodes A and C fulfils the condition. In addition, Table 11 exhibits the coordinates of the prescribed precisio positions for the synthesis in example 4.

For example 4, not only deformation energy shows a fast convergence but also the constraint error (see Fig. 17). In this graphic the values taken from Table 12 are plotted in linear scale as the behavior of the curves is clear enough. One can observe that although the number of iterations ascends to 9 because of the stop criterion used, acceptable convergence begins to settle from the 3rd iteration. Total time to reach the solution was 0.5 s in the two characteristics energy and error.

Table 13 displays the error produced in each synthesis prescribed point. One can observe the improvement achieved between initial mechanism and final one by the decrease of the error between position of the coupler E of the four-bar and the desired synthesis precision point. The data have been gathered from the kinematic analysis run with minimum distance function; in the one hand, with the initial mechanism of Fig. 15 and, in the other hand, with the final mechanism of Fig. 16.

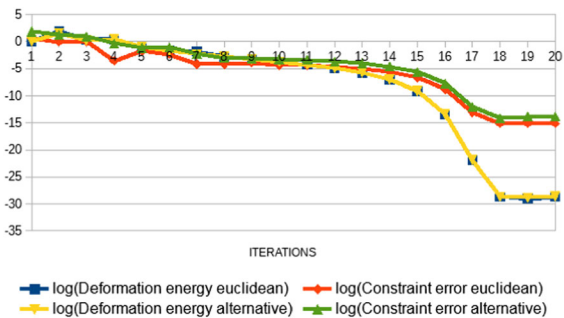

**Fig. 14** Evolution of the deformation energy, and of the error in the constraint with both approaches, example 3. The values have been displayed in logarithmic scale

The sequential movement of the mechanism of example 4 is shown in Fig. 18. The four-bar follows exactly the order of the prescribed points approximating the 11 desired positions quite accurately. The first precision point PP0 starts at the left end of the path, following them towards the right. It the case that the input link was required to be that marked with a 1 in Fig. 16, one would find that the resulting mechanism is in a near limit position in the last precision point. This could invalidate the usefulness of the mechanism, and currently no control of this kind of issues has been introduced in the algorithm. This is an interesting fact that is to be tackled in further developments.

Considering Table 13, the average distance in the initial mechanism is 0.03 units and in the final mechanism is 1.63e−05 units, much lower. With respect to the standard deviation, one can also state that for the initial guess it was 4.42e−02 units and for the final solution proposed it was improved to 1.34e−05.

**Table 8** Resulting values of deformation energy during the iterations for example 3 in both Euclidean and alternative approaches

| Iter | Euclidean fitness | log(Fitness) | alternative Fitness | log(Fitness) |
|------|-------------------|--------------|---------------------|--------------|
| 0    | 0.00              |              | 0.00                |              |
| 1    | 81.60             | 1.91         | 34.18               | 1.53         |
| 2    | 2.45              | 0.39         | 1.97                | 0.30         |
| 3    | 3.04              | 0.48         | 3.10                | 0.49         |
| 4    | 0.11              | −0.97        | 0.11                | −0.94        |
| 5    | 0.02              | −1.65        | 0.02                | −1.75        |
| 6    | 0.02              | −1.78        | 5.31E−03            | −2.27        |
| 7    | 1.77E−03          | −2.75        | 1.69E−03            | −2.77        |
| 8    | 5.72E−04          | −3.24        | 5.54E−04            | −3.26        |
| 9    | 1.78E−04          | −3.75        | 1.76E−04            | −3.75        |
| 10   | 5.22E−05          | −4.28        | 5.14E−05            | −4.29        |
| 11   | 1.23E−05          | −4.91        | 1.22E−05            | −4.91        |
| 12   | 1.89E−06          | −5.72        | 1.87E−06            | −5.73        |
| 13   | 1.11E−07          | −6.95        | 1.09E−07            | −6.96        |
| 14   | 7.47E−10          | −9.13        | 7.28E−10            | −9.14        |
| 15   | 4.35E−14          | −13.36       | 4.15E−14            | −13.38       |
| 16   | 1.52E−22          | −21.82       | 1.38E−22            | −21.86       |
| 17   | 2.13E−29          | −28.67       | 2.33E−29            | −28.63       |
| 18   | 8.48E−30          | −28.67       | 1.48E−29            | −28.83       |
| 19   | 2.03E−29          | −28.69       | 2.41E−29            | −28.62       |

**Table 9** Resulting values of constraint error during the iterations for example 3 in both Euclidean and alternative approaches

| Iter | Euclidean error | log(Error) | Alternative error | log(Error) |
|---|---|---|---|---|
| 0 | 4.72 | 0.67 | 69.39 | 1.84 |
| 1 | 0.99 | −4.35E−03 | 18.99 | 1.28 |
| 2 | 1.12 | 0.05 | 7.75 | 0.89 |
| 3 | 2.59E−04 | −3.59 | 0.47 | −0.33 |
| 4 | 0.02 | −1.76 | 0.08 | −1.12 |
| 5 | 4.47E−03 | −2.35 | 0.10 | −1.02 |
| 6 | 6.83E−05 | −4.17 | 5.17E−03 | −2.29 |
| 7 | 7.15E−05 | −4.15 | 9.21E−04 | −3.04 |
| 8 | 1.24E−04 | −3.91 | 6.89E−04 | −3.16 |
| 9 | 4.82E−05 | −4.32 | 4.89E−04 | −3.31 |
| 10 | 4.18E−05 | −4.38 | 3.88E−04 | −3.41 |
| 11 | 1.82E−05 | −4.74 | 1.94E−04 | −3.71 |
| 12 | 1.05E−05 | −4.98 | 9.92E−05 | −4.00 |
| 13 | 2.37E−06 | −5.63 | 2.45E−05 | −4.61 |
| 14 | 2.39E−07 | −6.62 | 2.30E−06 | −5.64 |
| 15 | 1.67E−09 | −8.78 | 1.66E−08 | −7.78 |
| 16 | 1.07E−13 | −12.97 | 9.95E−13 | −12.00 |
| 17 | 0.00 | | 7.11E−15 | −14.15 |
| 18 | 8.88E−16 | −15.05 | 0.00 | |
| 19 | 8.88E−16 | −15.05 | 1.42E−14 | −13.85 |

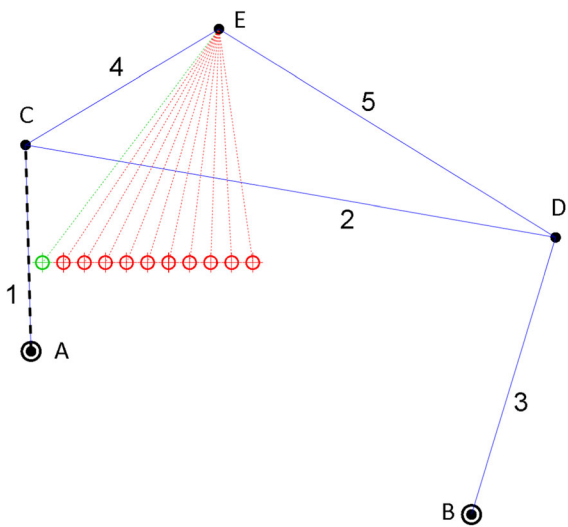

**Fig. 15** Initial mechanism with paths of prescribed synthesis points and constraint of fixed length, example 4

## 6.2 Topology with more links

In this example of a double butterfly shown in Fig. 19, it was aimed to draw a path in the shape of a hook for node K. This path was described by the 12 prescribed precision points defined in Table 14. During the path synthesis the element placed between nodes E and F was restricted to a fixed value of D=1.26 units.

**Table 10** Nodal coordinates of example 4 in initial and final configurations

| | $x_A$ | $y_A$ | $x_B$ | $y_B$ | $x_C$ | $y_C$ |
|---|---|---|---|---|---|---|
| Initial | −0.67 | −1.11 | 9.80 | −5.00 | −0.80 | 3.80 |
| Final | 1.53 | −4.49 | 8.32 | −6.36 | 1.14 | −0.51 |

| | $x_D$ | $y_D$ | $x_E$ | $y_E$ |
|---|---|---|---|---|
| Initial | 11.80 | 1.60 | 3.80 | 6.56 |
| Final | 12.35 | −2.31 | 7.00 | 7.05 |

Hence, the optimization problem for example 5 could be raised by Eq. (28) as:

$$\min_X f(\{X\}) = \sum_{j=1}^{12} \sum_{i=1}^{15} \left[ l_{ij}(\{x\}_j) - L_{ij}(\{X\}_j) \right]^2$$

subject to: $c_t(\{X\}) = L_t(\{X\}) - 1.26 = 0$

$$(28)$$

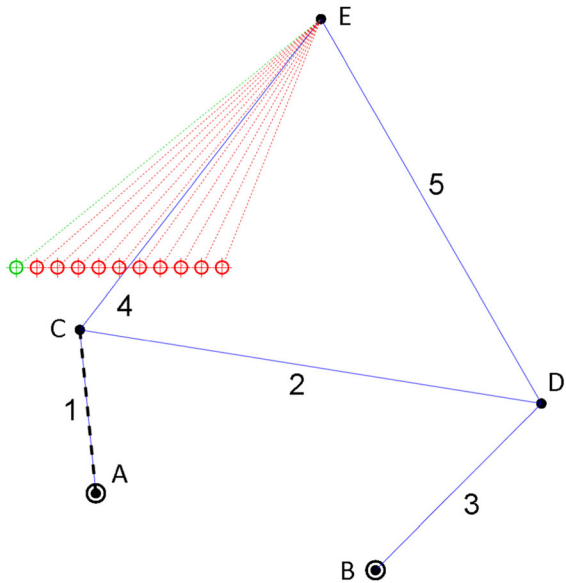

**Fig. 16** Solution mechanism, example 4

**Table 11** Coordinates of the prescribed precision points for the mechanism in example 4

|         | PP0  | PP1 | PP2 | PP3 | PP4 | PP5 | PP6 |
|---------|------|-----|-----|-----|-----|-----|-----|
| x coord | −0.4 | 0.1 | 0.6 | 1.1 | 1.6 | 2.1 | 2.6 |
| y coord | 1.0  | 1.0 | 1.0 | 1.0 | 1.0 | 1.0 | 1.0 |

|         | PP7 | PP8 | PP9 | PP10 |
|---------|-----|-----|-----|------|
| x coord | 3.1 | 3.6 | 4.1 | 4.6  |
| y coord | 1.0 | 1.0 | 1.0 | 1.0  |

(*) NOTE: PP = Precision Point

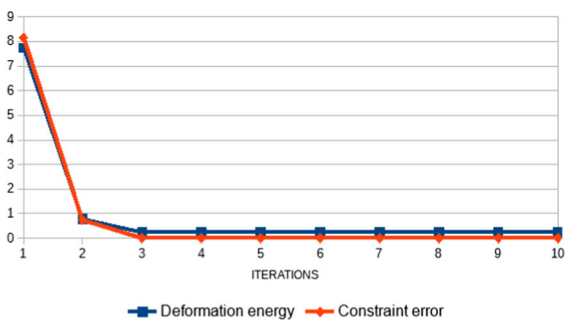

**Fig. 17** Evolution of the deformation energy (blue) and error in the constraint (red) with the Euclidean length approach, example 4. The values in this case have been displayed in linear scale. (Color figure online)

**Table 12** Resulting values of error and energy during the iterations for example 4 in the Euclidean approach

| Iter | Euclidean fitness | log(Fitness) | Error    | log(Error) |
|------|-------------------|--------------|----------|------------|
| 0    | 7.75              | 0.89         | 8.15     | 0.91       |
| 1    | 0.78              | −0.11        | 0.75     | −0.12      |
| 2    | 0.26              | −0.59        | 0.01     | −2.00      |
| 3    | 0.25              | −0.61        | 6.52E−05 | −4.19      |
| 4    | 0.25              | −0.61        | 2.3E−08  | −7.63      |
| 5    | 0.25              | −0.61        | 2.51E−06 | −5.60      |
| 6    | 0.24              | −0.61        | 2.46E−06 | −5.61      |
| 7    | 0.24              | −0.61        | 2.31E−08 | −7.34      |
| 8    | 0.24              | −0.61        | 3.55E−15 | −14.45     |
| 9    | 0.24              | −0.61        | 0.00     |            |

**Table 13** Errors committed at each point (distance between coupler point E and synthesis point) in example 4

|         | PP0      | PP1      | PP2      | PP3      |
|---------|----------|----------|----------|----------|
| Initial | 0.11     | 0.02     | 1.03e−04 | 1.49e−03 |
| Final   | 4.48e−05 | 1.57e−05 | 3.71e−05 | 1.35e−05 |

|         | PP4      | PP5      | PP6      | PP7      |
|---------|----------|----------|----------|----------|
| Initial | 1.38e−03 | 6.38e−06 | 2.86e−03 | 0.01     |
| Final   | 7.28e−08 | 1.28e−05 | 2.40e−05 | 1.37e−05 |

|         | PP8      | PP9      | PP10     |
|---------|----------|----------|----------|
| Initial | 0.04     | 0.08     | 0.12     |
| Final   | 2.19e−07 | 8.94e−06 | 8.06e−06 |

(*) NOTE: PP=Precision Point

The problem in exercise 5 can be described by the location of every prescribed precision point demonstrated in Table 14 and by the definition of all the nodal coordinates exposed in Table 15.

The coordinates shown in Table 15 only refer to the nodal positions, but with the help of Fig. 19 one can also understand the assembly configuration of the initial guess. The same can be stated for the solution configuration shown in Fig. 20.

Table 16 displays the results of some of the iterations necessary to achieve convergence. In total

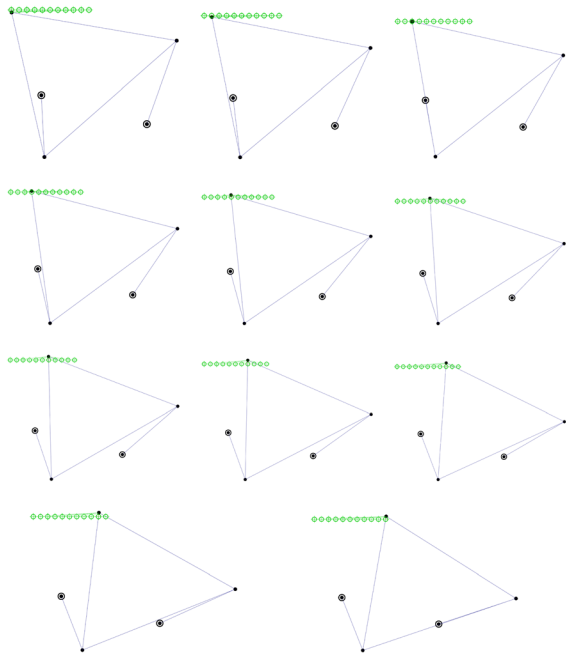

**Fig. 18** Fourbar of example 4 in the 11 prescribed synthesis points

**Table 14** Coordinates of the prescribed precision points for the mechanism in example 5

|        | PP0 | PP1 | PP2 | PP3 | PP4 | PP5 | PP6 |
|--------|-----|-----|-----|-----|-----|-----|-----|
| x coord | 4.6 | 4.4 | 4.1 | 3.6 | 3.1 | 2.6 | 2.4 |
| y coord | 7.6 | 8.2 | 8.6 | 8.9 | 9.0 | 8.6 | 8.2 |

|        | PP7 | PP8 | PP9 | PP10 | PP11 |
|--------|-----|-----|-----|------|------|
| x coord | 2.3 | 2.3 | 2.3 | 2.3  | 2.3  |
| y coord | 7.6 | 7.0 | 6.4 | 5.8  | 5.0  |

(*) NOTE: PP = Precision Point

**Table 15** Nodal coordinates of example 5 in initial and final configurations

|         | $x_A$  | $y_A$  | $x_B$ | $y_B$  | $x_C$ | $y_C$ |
|---------|--------|--------|-------|--------|-------|-------|
| Initial | −0.54  | −0.10  | 5.10  | −1.95  | 4.10  | 10.44 |
| Final   | 1.12   | −0.72  | 7.91  | −1.62  | 7.27  | 10.15 |

|         | $x_D$  | $y_D$ | $x_E$ | $y_E$ | $x_F$ | $y_F$ |
|---------|--------|-------|-------|-------|-------|-------|
| Initial | −0.37  | 2.06  | 1.44  | 2.36  | 1.62  | 3.62  |
| Final   | −3.70  | 2.20  | 6.24  | 1.64  | 5.47  | 2.64  |

|         | $x_G$ | $y_G$ | $x_H$ | $y_H$ | $x_I$ | $y_I$ |
|---------|-------|-------|-------|-------|-------|-------|
| Initial | 6.22  | 1.77  | 8.67  | 5.88  | 7.88  | 8.38  |
| Final   | 10.12 | 4.59  | 10.72 | 5.05  | 7.28  | 6.79  |

|         | $x_J$ | $y_J$ | $x_K$ | $y_K$ |
|---------|-------|-------|-------|-------|
| Initial | 6.84  | 5.90  | 4.08  | 8.02  |
| Final   | 9.32  | 6.02  | 4.60  | 7.61  |

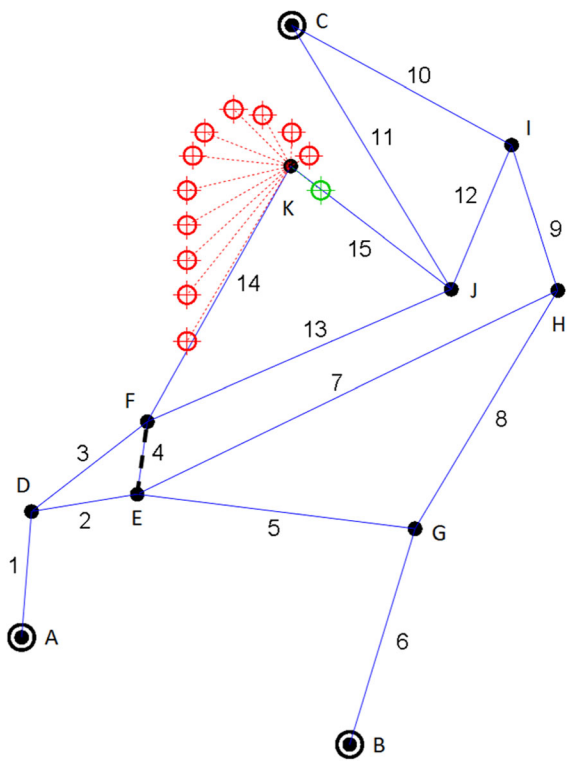

**Fig. 19** Initial mechanism with path of prescribed precision points and constraint of fixed length, example 5

there were needed 102 iterations, but near the 8th iteration in the case of the deformation energy, and around the 10th iteration in the case of the constraint error, reasonable convergence was achieved showing a fast convergence during the first steps. All these are shown in Fig. 21, where it can be concluded that the result was quite remarkable for example 5, not taking as many iterations as could be expected for such a complex mechanism with so many coordinates to relocate. In this case both the deformation energy and the error of the constraint are characterized by a behavior more similar to that in the previous examples. However, the criterion used in the algorithm to stop the process searches for such a high precision that the calculations continue until the 102 steps. The big amount of steps along with the complexity of the problems raised the time to solve it to about 200 s.

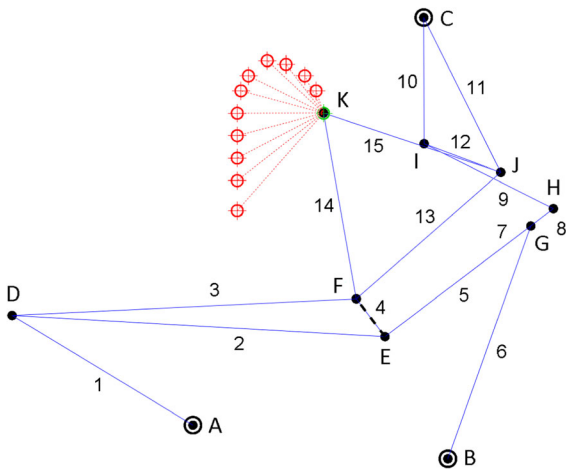

**Fig. 20** Solution mechanism, example 5

**Table 16** Resulting values of error and energy during the iterations for example 5 in the Euclidean approach

| Iter | Euclidean fitness | log(Fitness) | Error | log(Error) |
|---|---|---|---|---|
| 0 | 1.44 | 0.16 | 1.38E−05 | − 4.86 |
| 2 | 0.80 | − 0.10 | 0.02 | − 1.78 |
| 4 | 0.11 | − 0.94 | 0.65 | − 0.19 |
| 6 | 0.02 | − 1.82 | 2.01 | 0.30 |
| 8 | 1.56E−04 | − 3.81 | 0.01 | − 1.87 |
| 10 | 1.98E−05 | − 4.70 | 3.09E−04 | − 3.51 |
| 20 | 1.00E−08 | − 8.00 | 7.28E−04 | − 3.14 |
| 30 | 6.95E−10 | − 9.16 | 1.90E−05 | − 4.72 |
| 40 | 6.04E−10 | − 9.22 | 3.94E−07 | − 6.40 |
| 50 | 5.61E−10 | − 9.25 | 7.15E−07 | − 6.15 |
| 60 | 4.91E−10 | − 9.31 | 2.00E−06 | − 5.70 |
| 70 | 4.57E−10 | − 9.34 | 4.06E−06 | − 5.39 |
| 80 | 3.96E−10 | − 9.40 | 9.21E−07 | − 6.04 |
| 90 | 3.59E−10 | − 9.44 | 3.05E−07 | − 6.52 |
| 102 | 2.96E−10 | − 9.53 | 5.06E−07 | − 6.30 |

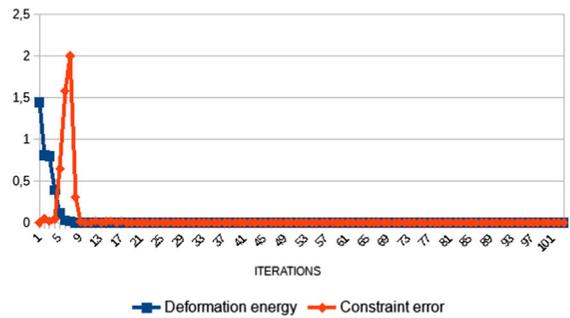

**Fig. 21** Evolution of the deformation energy (red) and error in the constraint (blue) with the Euclidean length approach, example 5. The values in this case have been displayed in linear scale

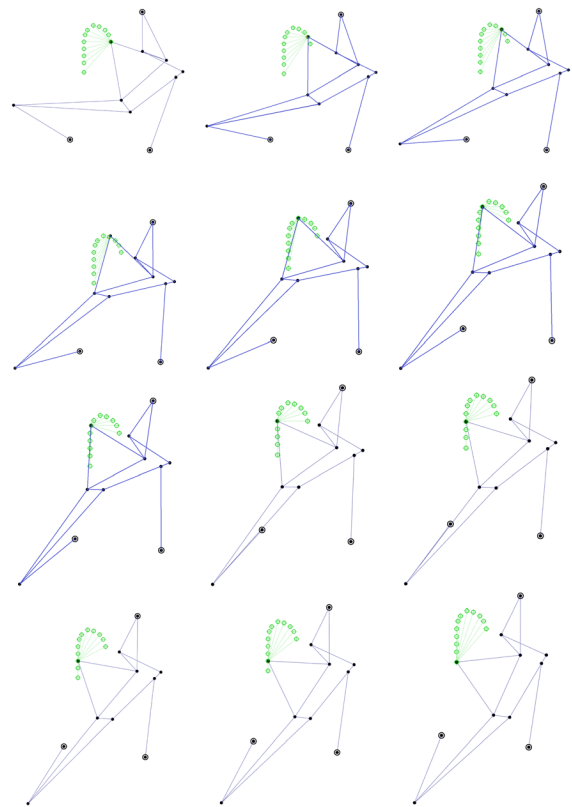

**Fig. 22** Double butterfly of example 5 in the 12 prescribed synthesis points

Figure 22 shows how the sequence followed by the movement of the mechanism corresponds to the order of the 12 prescribed points with a high accuracy.

With respect to the error produced in each synthesis point, one can observe the improvement achieved between initial mechanism and final one, as shown in

**Table 17** Errors committed at each point (distance between coupler and synthesis point), in example 5

|         | PP0       | PP1      | PP2      | PP3     |
| ------- | --------- | -------- | -------- | ------- |
| Initial | 8.03e−03  | 0.10     | 0.65     | 1.67    |
| Final   | 5.43e−08  | 6.22e−04 | 2.93e−05 | 0.02    |

|         | PP4      | PP5      | PP6      | PP7      |
| ------- | -------- | -------- | -------- | -------- |
| Initial | 2.73     | 2.22     | 2.72     | 2.94     |
| Final   | 8.96e−04 | 5.86e−04 | 1.39e−04 | 3.67e−05 |

|         | PP8      | PP9      | PP10     | PP11     |
| ------- | -------- | -------- | -------- | -------- |
| Initial | 1.82     | 1.62     | 1.38     | 1.29     |
| Final   | 8.74e−07 | 1.10e−04 | 1.00e−04 | 7.87e−10 |

(*) NOTE: PP = Precision Point

Table 17. The data have been gathered from the kinematic analysis with minimum distance function for the couple of mechanisms, the initial one in Fig. 19 and the final one in Fig. 20.

Therefore, the average distance in the initial mechanism is 1.59 units and in the final mechanism is 1.74e−03 units, much lower. Considering the standard deviation, it can be said that for the initial guess it was 9.40e−01 units and for the final solution proposed it was improved to 5.47e−03.

# 7 Conclusions

Two different types of features have been implemented to the general method of synthesis based upon the deformation energy. The first one allows the introduction of distance requirements between floating nodes. Two different approaches have been formulated along with their analytical derivatives. Both of them have shown similar performance. This has been verified by a set of examples developed ad-hoc. The Euclidean approach has a more direct physical interpretation as it is indeed the definition of distance. However, the derivatives in the alternative approach are easier to obtain as they do not have a square root of the Euclidean formula bringing a lower computational cost.

The second feature consist on the introducion of a constraint in the synthesis variables. This allows the definition of fixed lengths of links. Again, two approaches have been formulated. Both of them have shown to be effective, yielding satisfactory results as shown in the examples.

In order to introduce these features, nonlinear restrictions had to be introduced in the original optimization method, both for the analysis and the synthesis. This has been done using Lagrange multipliers, although other methods such as penalty functions or augmented Lagrange multipliers can easily be implemented. This allows to introduce more features in a near future with a small development cost.

**Acknowledgements** The authors thank the Spanish Ministry of Economy and Competitiveness for its support through Grants DPI2013-46329-P and DPI2016-80372-R (AEI/FEDER, UE), and the Education Department of the Basque Government for its support through Grant IT947-16.

**Compliance with ethical standards**

**Conflict of Interest** The authors declare that they have no conflict of interest.